\documentclass[12pt]{article}
\usepackage{amsmath}
\usepackage{amssymb}
\newcommand{\ds}{\displaystyle}

\begin{document}

\newtheorem{theorem}{Theorem}[section]
\newtheorem{remark}[theorem]{Remark}
\newtheorem{mtheorem}[theorem]{Main Theorem}
\newtheorem{observation}[theorem]{Observation}
\newtheorem{proposition}[theorem]{Proposition}
\newtheorem{lemma}[theorem]{Lemma}
\newtheorem{note}[theorem]{}
\newtheorem{extlemma}[theorem]{Ext-Lemma}
\newtheorem{corollary}[theorem]{Corollary}
\newtheorem{example}[theorem]{Example}
\newtheorem{definition}[theorem]{Definition}

\renewcommand{\labelenumi}{(\roman{enumi})}
\newcommand{\dach}[1]{\hat{\vphantom{#1}}}
\numberwithin{equation}{section}

\def\Z{{ \mathbb Z}}
\def\N{{ \mathbb N}}
\def\BB{ \mathbb B}
\def\DD{ \mathbb D}
\def\GG{ \mathbb G}
\def\BBB{ B_\BB}
\def\R{{\bf R}}
\def\D{{\hat{D}}}
\def\Q{{\mathbb Q}}
\def\G{\hat{G}}
\def\C{\hat{C}}
\def\T{{\cal T}}
\def\V{{\mathfrak V}}
\def\C{{\mathfrak C}}
\def\X{{\mathfrak X}}
\def\Y{{\mathfrak Y}}
\def\R{{\bf R}}
\def\D{\widehat{D}}
\def\A{\widehat{A}}
\def\G{\hat{G}}
\def\T{{\cal T}}
\def\B{\widehat{B}}
\def\BC{\widehat{B_C}}
\def\BCC{\widehat{B_\C}}
\def\restr{\restriction}
\def\Aut{{\rm Aut\,}}
\def\Im{{\rm Im\,}}
\def\ker{{\rm ker\,}}
\def\inf{{\rm inf\,}}
\def\sup{{\rm sup\,}}
\def\Br{{\rm Br\,}}
\def\Yphi{Y_{[\phi]}}
\def\Ypsi{Y_{[\psi]}}
\def\Xphi{X_{\tilde{\phi}}}
\def\Xpsi{X_{\tilde{\psi}}}
\def\a{\alpha}
\def\abar{\overline{\alpha}}
\def\aa{{\bf a}}
\def\to{\rightarrow}
\def\arr{\longrightarrow}
\def\sigmaa{{\bf \Sigma_a}}

\def\End{{\rm End\,}}
\def\Ines{{\rm Ines\,}}
\def\Hom{{\rm Hom\,}}

\def\restr{\upharpoonright}
\def\Ext{{\rm Ext}\,}
\def\Hom{{\rm Hom}\,}
\def\End{{\rm End}\,}
\def\Aut{{\rm Aut}\,}
\def\ker{{\rm ker}\,}
\def\defe{{\rm def}\,}
\def\rk{{\rm rk}\,}
\def\crk{{\rm crk}\,}
\def\nuc{{\rm nuc}\,}
\def\Dom{{\rm Dom}\,}
\def\Im{{\rm Im}\,}
\def\Yphi{Y_{[\phi]}}
\def\Ypsi{Y_{[\psi]}}
\def\Xphi{X_{\tilde{\phi}}}
\def\Xpsi{X_{\tilde{\psi}}}
\def\a{\alpha}
\def\abar{\overline{\alpha}}
\def\aa{{\bf a}}
\def\ra{\rightarrow}
\def\arr{\longrightarrow}
\def\iff{\Longleftrightarrow}
\def\sigmaa{{\bf \Sigma_a}}
\def\mm{{\mathfrak m}}
\def\F{{\mathfrak F}}
\def\X{{\mathfrak X}}
\def\Diam{\diamondsuit}
\def\mapdown#1{\Big\downarrow\rlap{$\vcenter{\hbox{$\scriptstyle#1$}}$}}
\def\cRp{\rm{cRep}$_2R$ }
\def\cR{c$R_2$}                  

\title{{\sc COTORSION THEORIES AND SPLITTERS}
\footnotetext{This work is supported by the project 
No. G-0294-081.06/93 of the German-Israeli
Foundation for Scientific Research \& Development\\
AMS subject classification:\\ 
primary 13D30, 18E40, 18G05, 20K20, 20K35, 20K40 \\
secondary: 13L03, 18G25, 20K25, 20K26, 20K30, 13C10 \\
Key words and phrases: cotorsion theories, completions, self-splitting modules, enough
projectives, realizing rings as endomorphism rings of self-splitting
modules\\
GbSh 647 in Shelah's list of publications}
} 

\author{ R\"udiger G\"obel and Saharon Shelah}

\date{}

\maketitle

\begin{abstract}
Let $R$ be a subring of the rationals. We want to investigate self
splitting $R$-modules $G$ that is $\Ext_R(G,G) = 0$ holds and follow Schultz 
\cite{Sch} to call such modules splitters. Free modules and torsion-free 
cotorsion modules are classical examples for splitters. Are there others? 
Answering an open problem by Schultz \cite{Sch} we will show that there are 
more splitters, in fact we are able to prescribe their endomorphism 
$R$-algebras with a free $R$-module structure. As a byproduct we are able to 
answer a problem of Salce \cite{Sa} showing that all rational cotorsion 
theories have enough injectives and enough projectives. 
\end{abstract}

\section{Introduction}

Jutta Hausen \cite{Ha} showed in her PhD-thesis under supervision of Reinhold
Baer in 1967 that any countable, torsion-free abelian 
group $G$ with $\Ext (G,G) = 0$
is free over some ring $R \subseteq \Q$. This interesting result, at that
time motivated by studies on automorphism groups, received new support
recently by investigations of cotorsion theories, see Salce \cite{Sa} and
Schultz \cite{Sch}. Our paper will deal with Hausen's result, that is with
groups $G$ such that $\Ext (G,G) = 0$. Following Schultz \cite{Sch},
we call these modules splitters. Because of their importance, also
other names are in use, see the  `Dictionary' on p. 351 in Ringel
\cite{Ri}. These self splitting modules are also called `stones' 
by Kerner \cite{Ke} (in contrast to     `bricks' which refers to the case
that the endomorphism ring is a division ring, see \cite{Ri1}), 
exceptional modules by Rudakov \cite{Ru} and Schur modules in 
Unger \cite{Un},  see also \cite{Wa}. We will stick to the name `splitters'
introduced by Schultz \cite{Sch}.
Our interest in splitters comes partly from their relatives, Whitehead 
groups, which
were investigated thoroughly over the last two decades, see results in
\cite{EM} and in a more recent paper \cite{BFS}. On the other
hand we are motivated in the study of splitters by open problems concerning
the algebraic structure of such modules and questions on cotorsion theories
related with splitters. 

\bigskip

First we notice a big difference between Whitehead
groups and splitters making splitters more attractive. In the case of
Whitehead groups one of the components in $\Ext(\_,\_)$ is 
$\Z$, hence countable!

\bigskip

Before we state our main results we will discuss the connection of cotorsion
theories and splitters, see $\S$ 2 for more details. Cotorsion theories were
introduced by Salce \cite{Sa} in 1978. They are dual to the
well-known classical torsion theories
replacing the crucial Hom-functor by $\Ext$. 
By reasons which will be clear soon, even if we only
want to know about splitters in the category of abelian groups 
we are bound to consider and will 
restrict to $R$-modules over subrings $R$ of the rationals $\Q$.

\bigskip

Following \cite{Sa}, a pair $(\F, \C)$ of maximal classes $\F, \C$ of 
$R$-modules is a cotorsion 
theory if $\Ext(\F,\C) = 0 $ in an obvious sense. Then
 $\C$ is the cotorsion
part and $\F$ is the torsion-free part of this theory $(\F, \C)$.
The most important example of a cotorsion theory is the `classical' cotorsion
 theory, developed in the 60th by Harrison and many
other algebraists, 
is the pair (torsion-free, cotorsion), see Fuchs \cite{Fu}.
It is easy verified that in the classical cotorsion theory the cotorsion
groups $C$ come from the rationals $\Q$ in the sense that
$$\Ext(C,\Q) = 0 \mbox{ if and only if } C \in \C,$$
so $\Q$ cogenerates the cotorsion theory.
 Cotorsion theories cogenerated by subgroups of $\Q$ are called
{\it rational cotorsion theories}; they are well-studied and the main 
objects in Salce's paper \cite{Sa}. From the homological point of view it is 
important to know, whether rational cotorsion theories have enough 
projectives and injectives. Only then we are ready to introduce cotorsion 
hulls! This question was raised in Salce \cite{Sa} and remained open so far. 
We will answer it positively in Section 6. 

\medskip
The answer to Salce's question is a byproduct of our study of splitters. 
As already indicated, the study of abelian group splitters
can be reduced to torsion-free, reduced $R$-modules $(R \subseteq \Q)$ by a
``Reduction Theorem \ref{2.6}'' due to Schultz \cite{Sch}, where $R = \nuc G$ is the
nucleus of $G$. The nucleus of $G$ is the largest subring $R$ of $\Q$ such
that $G$ is (canonically) an $R$-module. The first step towards the structure
of splitters is Hausen's theorem which can be slightly extended (Theorem \ref{2.7}):

\medskip
{\it Splitters of cardinality $< 2^{\aleph_0}$ are $\aleph_1$-free modules
over their nuclei}.

\medskip
The key tool of this paper can be found in
Section 3. Here we will prove our Ext-Lemma \ref{3.4}.
Our original proof of the Ext-Lemma was based on constructing solutions of
of certain systems of equations. However, following 
a successor paper \cite{GT} we will 
replace this by a homological shortcut. Nevertheless
homological arguments often conceal details of the  
structure under investigation even though the arguments may
faster lead to the desired result. The Whitehead
problem only settled after two decades may explain this very well; this is
one reason why we also include a non-homological proof (\ref{4.10}) of 
(\ref{2.7}) as well. This will  also shed light on the Ext-Lemma \ref{3.4}. 

\medskip

Obviously free groups and torsion-free cotorsion groups are
splitters, the first by trivial algebraic reasons the latter by completion.
Before we now state the Ext-Lemma 
we also need a definition of $n$-free-by-1 $R$-modules and we explain its
connection with classical results. Similar to simply presented groups
$n$-free-by-1 groups are easy represented by free generators and
relations. 
If we want to find non-free splitters $G$, which are necessarily of 
cardinality $\geq \aleph_1$ by (\ref{2.7}), one of the obstacles 
are ``small'' non-free subgroups of $G$. If $G$ has non-free countable 
subgroups, then by Pontryagin's theorem
$G$ also has non-free subgroups of some minimal finite rank $n+1$. 
These groups are investigated under the name $n$-free-by-1 $R$-modules in 
Section 3 first. The name is easily explained; see Observation \ref{3.1} 
and (3.2) for their elementary
properties. Rational groups are the special 0-free-by-1 $\Z$-modules, but
those for $n \geq 1$ or $n = \omega$ are particular important. Such groups
$G'$ are canonically connected with certain ``easy'' systems of 
equations (\ref{3.3}) and another group $G$ is called $G'$-complete if these 
equations are always solvable in $G$. This observation surely is connected 
with ideas of type in model theory, see Prest \cite{Pr}. 
Our crucial $\Ext$-Lemma now reads as follows:

\bigskip
{\bf Ext-Lemma \ref{3.4}} {\it Let $G'$ be an $n$-free-by-1 $R$-module. Any 
torsion-free module $G$ over the nucleus $R$ is $G'$-complete if and only if
$\Ext (G',G) = 0$}. 

\medskip

As indicated, $\Z$-adically complete modules are $G'$-complete modules 
for a suitable module $G'$.
If $\Ext (\Q, G) = 0$ and $G$ is torsion-free, reduced,
then $G$ is complete in the $\Z$-adic topology by classical results due to
Kaplansky, see \cite{Fu}. Hence splitting quite often implies completion.
Such torsion-free modules are also cotorsion
modules, hence algebraically compact. They have a nice structure
theory by cardinal invariants studied by many algebraists, see Warfield
\cite{War}, Ziegler \cite{Zi} and work by {\L}o\v{s}, Kaplansky and others, 
see e.g. \cite{Fu, EM}. Moreover, they are the cotorsion-splitters of 
the classical cotorsion theory. 
If $\Q$ is replaced by another rational group, this can also be seen in Salce
\cite{Sa}. 
Note that the structure theory for these rational cotorsion theories
by Salce \cite{Sa} extends easily for $n$-free-by-1 $R$-modules.

\medskip

In the first part of Section 4 we deduce classical results,
like Kaplansky's theorem (above) from the $\Ext$-Lemma. Then we clarify the
new situation shown by (\ref{3.4}): 
$G'$-complete $R$-modules do not need to be
complete in any $S$-adic topology.
There are $n$-free counterexamples (where any submodule of rank $< n$ is 
free).
 Their endomorphism ring (with free additive
structure) may be prescribed as
well (Theorem \ref{4.7}). This shows that in contrast to classical 
$\Z$-adically complete groups these $G'$-complete $R$-modules $G$ cannot
be classified by any reasonable invariants. This is the case even if $G$ has
a nice-looking filtration build up by copies of $G'$ (Theorem \ref{4.9}). 
Some of
these $n$-free examples turn out to be splitters as shown in Section 5 
and we are able 
to prescribe endomorphism rings of splitters. 
We will see that for a
given number $n >1$ there are $n$-free splitters
of size $2^{\aleph_0}$ (and larger), which answers an open problem.
These $n$-free modules are obviously cotorsion-free.
Modifying our arguments, Eklof noticed that $n$-free can be replaced by
$\aleph_1-$free if the size of the splitter is at least
$2^{\aleph_1}$. Changing the construction again slightly as in
\cite{GT} it is now also possible to show that there are
non-free but $\aleph_1$-free and slender splitters of this cardinality, see
\cite{GT}. Hence it should be interesting to study $\aleph_1$-free
splitters of cardinality $\aleph_1$. This earlier part of Section 7 grew to
an individual joint paper \cite{GS}. We would like to thank Paul Eklof for
pointing out a wrong argument in that first version part of Section 7,
which as a consequence gave rise to \cite{GS}. Here we show that these
modules are free indeed in ordinary set theory (ZFC).

\bigskip

\section{Cotorsion theories and splitters - a summary of known facts 
and definitions; some generalizations}

In 1966 S. E. Dickson \cite{D} introduced torsion theories for abelian 
categories by exploiting the $\Hom$-functor. This helped to overcome
difficulties in defining torsion submodules and provided grounds for
new research.    Replacing formally the
Hom-functor by the Ext-functor, Salce \cite{Sa} developed the basic
tools for a cotorsion theory, which naturally extends the `classical
cotorsion theory', the notion of cotorsion (modules) abelian groups 
introduced by Harrison,  Nunke,  Fuchs,  Kaplansky and others,
see Fuchs \cite[p.232]{Fu}. We will use some notation from G\"obel, Prelle
\cite{GP} to introduce Salce's cotorsion-theory. 
If $\X, \Y$ are two classes of abelian groups, we say that\\
$(i)$ \quad $\X \perp \Y$ if and only if $\Ext (X,Y) = 0$ for all 
$X \in \X$ and $Y \in \Y$.\\
Moreover\\
$(ii)$ \quad $\X^\perp$ is the unique largest class of abelian groups with
$\X \perp \X^\perp$ \\
and dually \\
$(iii)$ \quad$^\perp \X$ is the unique largest class of abelian groups with
$^\perp \X \perp \X$.\\ 
The class $\X^\perp$ is called the injective
closure of $\X$ and $^\perp \X$ is the projective closure of
$\X$, notions which are natural in view of Salce \cite{Sa}. Following
Salce \cite{Sa} we say that a pair $(\F, \C)$ of classes $\F,
\C$ of abelian groups is a cotorsion theory if the following conditions
hold.\\
$(iv)$ \quad $ 
\begin{cases} (a) \quad \F \perp \C \\
(b) \quad \mbox{The classes are maximal: } \  \C = \F^\perp
\ \mbox{ and }  \F = {}^\perp\C.
\end{cases} $

\noindent 
$\C$ is the cotorsion-part and $ \F $  is the torsion-free 
class of this theory. The `classical cotorsion theory' 
$(\F_c, \C_c)$ deals with
the torsion-free part $\F_c =  \mbox{ all torsion-free groups, }$
and the cotorsion-part $ \C_c = \mbox{ all cotorsion groups}$.

There are two dual ways to produce new cotorsion theories from a given class
$\X$ of abelian groups. We either begin with $^\perp \X$ or with
$\X^\perp$, respectively and get:

\bigskip
\noindent
$(v)$ \quad 
$\begin{cases} \X \;\ \mbox{{\bf generates} the cotorsion theory}\;\;
(^\perp \X, (^\perp \X)^\perp).\\
\mbox{Dually,}\;\ \X \;\ \mbox{{\bf cogenerates} the cotorsion
theory}\;\ (^\perp (\X^\perp), \X^\perp).
\end{cases}$

\bigskip
It is easy to check that the pairs in $(v)$ are cotorsion theories
and the `classical cotorsion theory' is cogenerated by the rationals $\Q$, 
i.e. $(\F_c, \C_c) = (^\perp(\Q^\perp), \Q^\perp)$. The crucial part is
$\C_c = \Q^\perp$ i.e. $\Ext (\Q,G) = 0$ for all cotorsion modules
$G$. This well-known fact \cite{Fu} is basic in $\S$3. The well studied
classical cotorsion theory suggests a close investigation of their immediate
relatives coming from subgroups of $\Q$: \\
Suppose $S \subseteq \Q$ is a
$rank$-1 group and assume that $1 \in S$ without loss of generality. The main
task of Salce \cite{Sa} is a detailed description of the cotorsion theory
$(\F_S, \C_S)$ cogenerated by $S$. These cotorsion theories (for any 
$S \subseteq \Q)$ are called rational cotorsion theories. We will make 
use of Salce's main theorem describing the class $S^\perp$ of 
$S$-cotorsion groups in $(^\perp(S^\perp), S^\perp)$.
The following three standard definitions are needed:\\
$(vi)$ \quad $\chi_S = \chi (1) = (r_p)_{p \in \Pi}$ \\
(where $\Pi$ is the set of primes)
is the characteristic of $S$ with $p^{r_p}$ the maximal $p$-power dividing
1 in $S$, or $r_p = \infty$ if $p$ divides 1 infinitely often.\\
$(vii)$\quad $G_S = \bigcap\limits_{p \in \Pi} p^{r_p} G$\\
where $p^\infty G = \bigcap\limits_\sigma p^\sigma G$ is defined by the Ulm
sequence inductively for all ordinals:

$p^{\sigma+1} G = p (p^\sigma G)$ and $p^\sigma G = \bigcap\limits_{\nu < \sigma}
p^\nu G$ for limit ordinals $\sigma$.

\medskip
\noindent
$(viii)$ \quad If $r_p < \infty$, then $G_p^S = G/p^{r_p} G$ and $G_p^S =
\stackrel{\bullet}{G}_p =
\Ext (Z (p^\infty), G)$ if $r_p = \infty$.

\bigskip
\begin{theorem}[\cite{Sa}] \label{2.1}
If $S \subseteq \Q$ is as above, then \\
$G$ is $S$-cotorsion $\iff$ $G/G_S$ is cotorsion $\iff$
$G/G_S \cong \prod\limits_{p \in \Pi} G^S_p.$
\end{theorem}

\bigskip

For homological consideration it is very important to find out whether some
category has enough injectives and projectives. A cotorsion theory $(\F,
\C)$ has {\bf enough injectives} if and only if for all abelian groups 
$G$ there are $C \in \C, F \in \F$ and a short exact sequence
$$0 \ra G \ra C \ra F \ra 0.$$
Dually $(\F, \C)$ has {\bf enough projectives} if and only if 
there is another short exact sequence
$$0 \ra C \ra F \ra G \ra 0.$$
An easy lemma reduces the question about the existence of injectives 
and projectives to one problem.

\medskip
\begin{lemma}[\cite{Sa}] \label{2.2} A cotorsion theory $(\F, \C)$
has enough projectives if for all free groups $A$ there are $C \in \C,
F \in \F$ and a short exact sequence
$$0 \ra A \ra C \ra F \ra 0.$$
\end{lemma}
Hence $(\F, \C)$ has enough projectives if and only if it has enough
injectives.\\
Classical results, widely used, show that any abelian group $A$ can be purely
embedded into its cotorsion hull $A^\bullet$, see \cite[p. 248]{Fu}. Hence the
classical cotorsion theory has enough injectives and by the above enough
projectives. Naturally, Salce \cite[Problem 2, p. 31]{Sa} raised the question
whether rational cotorsion theories have enough projectives (injectives). We
will answer this question in the affirmative in a more general context - as
a by-product of our study of splitters.

\medskip
`Splitters' were introduced in Schultz \cite{Sch}. 
They also come up
under different names; see the introduction. Moreover we will see
immediately that splitters are closely connected with
Salce's work \cite{Sa}. Recall that the cotorsion class $\C$ of a
cotorsion theory $(\F, \C)$ is closed under epimorphic images. 
Similarly $\F$ is closed under subgroups. If $F \in \F, C \in \C$
are the groups in the exact sequence to define enough projectives, respectively
injectives for $G$, then either $F \in \F \cap \C$ or $C \in \F
\cap \C$. Hence $\F \cap \C$ is particular important. In the
classical case this is the class of torsion-free cotorsion groups or
equivalently torsion-free algebraically compact groups, which can be
classified by cardinal invariants, an extension to rational
cotorsion-groups is given in Salce \cite{Sa}.
Elements $G$ in $\F \cap \C$ obviously satisfy the condition
$\Ext (G, G) = 0$, i.e. the sequence\\
$(ix)$ \quad $0 \ra G \ra * \ra G \ra 0$ always splits 
and $G$ is self splitting or a splitter.

We arrive at the 

\medskip
\begin{definition}[\cite{Sch}] \label{2.4}
An $R$-module $G$ over a (hereditary)
commutative ring $R$ is an $R$-splitter if any $R$-module sequence
$(ix)$ splits
or equivalently $\Ext_R (G,G) = 0.$ 
\end{definition}

\medskip
We are mainly interested in subrings $R$ of $\Q$. In this case, if $G$ is a
torsion-free $R$-module, then 

$(x)$ $\Ext_R (G,G) = \Ext_\Z (G,G)$

\noindent
because $\Z$-homomorphisms are $R$-homomorphisms of $G$ and we call $G$ a 
splitter if (\ref{2.4}) holds. 
Obvious examples of splitters are the torsion-free
cotorsion groups in $\F_c \cap \C_c$, coming from the 
classical cotorsion theory $(\F_c,\C_c)$. The
other example comes from a trivial cotorsion theory 
$$(\F = \mbox{ free groups, } \C = \mbox{ all groups),} $$
hence free groups in $\F \cap \C$ are
splitters. In view of the countable case of (\ref{2.7}) and the above, Schultz
\cite[Problem 4]{Sch} raised the question whether these are all splitters.
We will answer this question to the negative in Section 5. 
However, following Schultz \cite{Sch}, we first reduce the problem to the
torsion-free case.

\bigskip
In order to investigate splitters, Schultz \cite{Sch} introduced the very
useful notion of a nucleus of a group.
\begin{definition} \label{2.5} The nucleus of a torsion-free group
$G \neq 0$ is the largest subring $R = \nuc G$ of $\Q$ such that $G$ is an
$R$-module. Hence $R$ is generated as a subring of $\Q$ by all $\frac{1}{p}$
($p$ any prime) for which $G$ is $p$-divisible, i.e. $p G = G$.
\end{definition}

\bigskip
We will fix this notion $R = \nuc G$ of a nucleus of $G$ throughout this paper.
The following result reduces the study of splitters among abelian
groups to those which are torsion-free and reduced modules over their nuclei.

\bigskip
\begin{theorem}[\cite{Sch}] \label{2.6}
Let $G$ be any abelian
group and write $G = D \oplus C$ as a decomposition of $G$ into the
maximal divisible subgroup $D$ and a reduced complement $C$.
Moreover let $\pi$ be the set of all those primes for which $D$ has a
non-trivial primary component. Then the following conditions are
equivalent.

\bigskip
$(i)$ $G$ is a splitter.

\bigskip
$(ii)$ $\begin{cases} 

\mbox{(a)}\;\  C \;\ \mbox{is a torsion-free (reduced) splitter with }
pC = C \mbox{ for all } p \in \pi. \\
\mbox{(b)} \mbox{ If } D \mbox{ is not torsion then } C
\mbox{ is cotorsion.}
\end{cases}$
\end{theorem}

Condition $(ii)(a)$ of the theorem say that $G$ is and $R$-module over
the ring generated by all primes $p^{-1}$ with $ \ p \in \pi$.
For convenience of the reader we sketch the essential steps of the proof.
(\ref{2.6}) is based on an easy observation, see \cite{Fu}:

\medskip
{\bf ($*$)} {\it If $\Ext (A,B) = 0, A' \subseteq A, B' \subseteq B$ then 
$\Ext (A', B/B') = 0$ as well.} 

\medskip
A short exact sequence
$$ 0 \arr B \overset{\beta}{\arr} C \overset{\alpha}{\arr} A \ra 0$$
represents $0$ in $\Ext(A,B)$ if and only if there is a splitting map
$\gamma: A \arr C$ such that $\gamma \alpha = id_A$. 

\bigskip
$(ii)$ $\ra$ $(i)$ is obvious and $(i)$ and ($*$) imply $\Ext (C,C)=0$. 
To prove that $C$ is torsion-free, assume to the contrary, that $C$
contains a copy $Z_p$ of a cyclic group of order $p$ for some prime 
$p \ne 1$ and consider two cases $pC = C$ and $pC \ne C$. The case
 $p C = C$ is impossible because $Z (p^\infty)$ would be a subgroup of $C$ 
but $C$ is reduced. If
$p C \neq C$ we find $C' \subset C$ with $C/C' \cong Z_p $. 
Observation ($*$) gives $0 = \Ext (Z_p, C/C') = \Ext (Z_p, Z_p) \cong Z_p$ a
contradiction; so $C$ is torsion-free.
If $D_p \neq 0$, then $\Ext (Z(p^\infty), C) = 0$ by ($*$),
hence $p C = C$ by \cite[Theorem 52.3]{Fu}. If $\Q \subseteq D$, then
$\Ext (\Q, C) = 0$ by ($*$), and $C$ must be cotorsion.

\medskip
From now on we will assume that $G$ is a torsion-free, reduced $R$-module
where $\nuc G = R$. A classical result due to Hausen \cite{Ha}
states that countable splitters $G$ are free $R$-modules indeed. This can be
slightly extended to say

\medskip
\begin{theorem} \label{2.7} If $R = \nuc G$ is the nucleus of the torsion-free
group $G$ and $G$ is a splitter of cardinality $< 2^{\aleph_0}$, then $G$ is
an $\aleph_1$-free $R$-module.
\end{theorem}

\medskip
Recall that $G$ is an $\aleph_1$-free $R$-module, if any countably generated
$R$-submodule is free. 
We will say that $G$ is $n$-free for some natural number $n$ if all the 
submodules of $G$ of rank $\le n$ are free. This agrees with our notion of 
$n-$free-by$-1$ $R$-modules in Section 3.
Theorem \ref{2.7} will be important in Section 7. We will
provide a quite obvious homological proof, followed later by some direct 
arguments leading to the same result. We believe that the direct arguments
uncover what's hidden by homology! The non-homological proof is at the end
of Section 4 in Corollary \ref{4.10}.

\medskip
{\bf First proof of (\ref{2.7})} (the homological approach). 
The proof follows in two
steps. It is convenient to recall ($*$). 
Also we say that an $R$-submodule $E$
of an $R$-module $C$ has {\bf full rank} if $C/E$ is torsion. First we claim

\medskip
(a) $\left\{
\begin{tabular}{p{12cm}}
{\it If $E$ is a full rank $R$-submodule of an $R$-submodule $C$  
of finite rank, both $E$ and $G$ have the same nucleus $R$, and $\Ext (C,G) = 0$,
then there is another full rank submodule $F$ with $E \subseteq F \subseteq C,
F/E$ finite and $\Ext (C/F,G) = 0$. Moreover, if $E$ is a free $R$-module, then
$F$ is free as well}.
\end{tabular}\right.$

\medskip
{\bf Proof.} We take dual groups $X^* = \Hom (X,G)$ and let
$$0 \ra E \ra C \ra C/E \ra 0$$
be the obvious short exact sequence. Hence
$$E^* \ra \Ext (C/E,G) \ra \Ext (C,G) = 0.$$
\medskip
However $|E^*| < 2^{\aleph_0}$ and $|\Ext (C/E,G)| < 2^{\aleph_0}$ follows.
The number of primes $p$ with $\Ext ((C/E)_p, G) \neq 0$ must be finite and 
the $p$-primary components $(C/E)_p$ of the torsion module $C/E$ must be
finite as well. Take $E \subseteq F \subseteq C$ with $F/E = \bigoplus\limits_p
(C/E)_p$ which is finite. We have $C/E = C/F \oplus F/E$ and $\Ext (C/E,G) =
\Ext (F/E,G)$, hence $\Ext (C/F,G) = 0$.

\medskip
If $E$ is free and $|F/E| = n$, then $F \cong n F \subseteq E$ is free.

\medskip
(b) $\left\{ \begin{tabular}{p{12cm}}
{\it If $\Ext (H,G) = 0$ and $H,G$ have the same nucleus $R$, then $H$ is
an $\aleph_1$-free $R$-module}.
\end{tabular}\right.$

\medskip
{\bf Proof.} Let $C \subseteq H$ be any $R$-submodule of finite rank. Hence
$\Ext (C,G) = 0$ by ($*$). Clearly there is a free $R$-submodule $E 
\subseteq C$
of full rank. From (a) we have $E \subseteq F \subseteq C$ with 
$F/E$ finite, $F$ free and $\Ext_R (C/F,G) = 0$. If there is $p$ with 
$(C/F)_p \neq 0$, 
then $p G = G$ is a contradiction for $R$, hence $C = F$ is free. 
Pontryagin's theorem completes the proof. $\hfill{\square}$

\section{The Ext-Lemma}

The main result of this section is related to a well-known 
observation due to Harrison and Kaplansky, see Fuchs \cite[p. 247--249]{Fu}. 
If $G$ is torsion-free and reduced, then $G$ is cotorsion if 
$\Ext (\Q,G) = 0$, and this is the same as
to say that $G$ is complete in the (Hausdorff) $\Z$-adic 
topology. Recall that
$G$ is cotorsion if and only if $\Ext (G',G) = 0$ for {\bf any} 
torsion-free group $G'$,
which explains the strength of the demand $\Ext (\Q,G) = 0$. How much of the
completion is left over, if $\Q$ is replaced by particular groups $G'$? If
$G'$ is a torsion-free group of rank 1, that is a subgroup of $\Q$, 
this question is answered by
Theorem 3.5  in Salce \cite[p.21]{Sa}, which is basic for his rational
cotorsion theories, see Section 4. Here we are interested in relatives of
these rank-1 groups, which occur naturally as subgroups of torsion-free 
groups. We begin with an easy motivation concerning these relatives of 
rank-1 groups by showing their existence as subgroups of arbitrary 
torsion-free groups.\\
\begin{observation} \label{3.1}
(a) If $G$ is torsion-free but not $\aleph_1$-free over its nucleus $R$,
then there is a pure $R$-submodule $G' \subseteq G$ of minimal finite rank,
which is not $R$-free.\\
(b) If $rk G' = n+1$ then we can find a free $R$-module 
$F = \bigoplus\limits_{m \in \omega} y''_m R \oplus 
\bigoplus\limits_{i<n} x_i''R $ and elements $k_{i m}, p_m \in R$ 
(the $p_m$'s constitute a divisibility chain: $p_j d_{jm} = p_m$ 
($j \le m$) for some $d_{jm}\in R $) with 
$$  G' \cong F/N \mbox{ if }
N = \left < y''_{m+1} p_m - y''_m - \sum\limits_{i<n} x''_i k_{i m}: 
m \in \omega \right >_R \subseteq F.$$
(c) $N = \bigoplus\limits_{m < n}(y''_{m+1} p_m - y''_m - 
\sum\limits_{i<n} x''_i k_{i m})R $ is a free $R$-module.
\end{observation}
{\bf Proof.}
(a) The nucleus $R$ is a PID, hence Pontryagin's theorem applies,
see Fuchs \cite{Fu}. There is a pure $R$-submodule $G'$ of finite rank which 
is not free. Clearly we may choose $G'$ of minimal rank $n+1$.

\medskip

(b) Let $G' = \left < x_0, \cdots, x_n\right >_*$ be the $R$-module of rank 
$n+1$ given by (a). By the minimality of $n$ we observe that 
$$\left < x_0, \cdots, x_{n-1} \right >_R = \bigoplus\limits_{i<n} x_i R$$ 
is a free, pure submodule of $G'$. If $G'$ has
rank 1, then this direct sum is zero. Also note that the 
torsion-free rank-1 module $G'/\bigoplus\limits_{i<n} x_i R \subseteq \Q$ is 
generated by $y_m \in G'$$(m \in \omega)$ and relations $y_0 = x_n$ and 
$y_{m+1} p_m \equiv y_m \mod\bigoplus\limits_{i<n} x_i R$, see Fuchs 
\cite[Vol. 2]{Fu}.\\Hence $y_{m+1} p_m - y_m \in \bigoplus\limits_{i<n} x_i 
R$, and there are $k_{i m} \in R$ $(i<n)$ such that
$$y_{m+1} p_m = y_m + \sum\limits_{i<n} x_i k_{i m}.$$
If $F$ is the free $R$-module given in the Observation \ref{3.1} (b), then
$$y''_m \ra y_m \;\; ,\;\;  x''_i \ra x_i \;\; (m \in \omega, i<n)$$
is a well-defined epimorphism from $F$ onto $G'$ with Kernel $N$ as in 
(\ref{3.1})(b).

(c) [Note that the proof remains valid if $n$ is replaced by 
$\omega$.] 
Consider elements $s_m \in R$ such that
$$\sum^k\limits_{m=0} (y''_{m+1} p_m - y''_m - 
\sum\limits_{i < n} x''_i k_{i m}) s_m = 0.$$ 
The coefficient of $y''_{k+1}$ is $p_k s_k = 0$, hence $s_k = 0$. 
Inductively we get $s_0 = \cdots = s_k = 0$ and the sum in (c) is direct.
$\hfill{\square}$

\bigskip
Observation \ref{3.1} explains our interest in torsion-free
$R$-modules of rank $n+1$, which are extensions of free $R$-modules 
of rank $n$ by a torsion-free $R$-module of rank 1. We will 
always assume that such a module is not free or equivalently that 
the rank-1 group is not $R$. 

We are also interested in rank-1 extensions of free $R$-modules of 
countable rank. More generally, let $n \le \omega$ and define 
$$B = \bigoplus\limits_{i < n} x_i R \oplus 
\bigoplus\limits_{m \in \omega} y_m R$$ 
be a free $R$-module. 
Consider elements 
$ p_m, k_{im}, d_{jm}$ as in Observation \ref{3.1} 
with $k_{im} = 0 $ for almost all $i$ and each $m$. Let
$$N = \left <   y_{m+1}p_m - y_m - \sum\limits_{i < n} x_i 
k_{i m}:\ \ m \in \omega  \right >_R$$
be a free submodule of $B$. The quotient module $G'= B/N$
 satisfies the relations
$$y'_{m+1} p_m = y'_m + \sum\limits_{i < n} x'_i k_{i m}, 
\ \  \  p_j d_{j m} = p_m \ \ (j \le m \in \omega),$$
where $y'_m = y_m + N$ and $x'_i = x_i + N$.

We will use these particular almost free $R$-modules, their representation 
and related `closures' very often and therefore summarize

\medskip
\begin{definition} \label{3.3}
 If $n \leq \omega$ and $R$ is a subring of $\Q$, then using $B$ and $N$ 
above, we define an $n$-free-by-1 $R$-modules $G'$ as the quotient module
$B/N$ or equivalently the module freely generated by
$$G' = \left < y'_m,\ \bigoplus\limits_{i<n} x'_i R: \quad m \in \omega 
\right >_R$$
except the relation 
$$y'_{m+1} p_m = y'_m + \sum\limits_{i<n} x'_i k_{i m} \quad (m \in \omega)$$
with $p_m$'s a divisibility chain as above.\\
Moreover, if $G$ is any torsion-free $R$-module over its
nucleus $R$, then we say that $G$ is $G'$-complete 
if for any sequence $c_m \in G$ $(m \in \omega)$ the system of 
equations
$$y_{m+1} p_m = y_m + \sum\limits_{i < n} x_i k_{im} + c_m \quad 
(m \in \omega)$$
has solutions $y_{m}, x_i \in G$ $(m \in \omega, i < n)$.
\end{definition}

As we assume that the divisibility chain of $p_m$'s in Definition \ref{3.3}
defines a proper type (not $R$), $0$-free-by-1 $R$-modules are 
the old-fashioned non-free, torsion-free rank-1 $R$-modules.
In Corollary \ref{4.1} we will show that $\Q$-completions are the well-known 
$\Z$-adic completions. The key for this paper is the following connection
between $\Ext$ and completions.

\begin{extlemma} \label{3.4}
Let $G$ be a torsion-free $R$-module over its nucleus $R$ and $G'$ be an 
$n$-free-by-$1$ $R$-module for some $n \leq \omega$. Then $\Ext (G', G) = 0$
if and only if $G$ is $G'$-complete. 
\end{extlemma}

{\bf Proof}. Let $n \leq \omega$ and let
$G' = \left < y'_m, \ \bigoplus\limits_{i < n}
x'_i R: \ \ m \in \omega \right >_R$
be expressed as in Definition \ref{3.3}.
We have a short exact sequence
$$  0 \arr N \overset{\sigma}{\arr} B \arr G' \arr 0$$
and it follows 
$$\Hom(G',G) \arr \Hom(B,G) \overset{\sigma^*}{\arr} \Hom(N,G) \arr
\Ext(G',G) \arr \Ext(B,G)$$
where $\Ext(B,G) = 0$ from freeness of $B$. Hence
$$\Ext (G',G) = 0 \mbox{ if and only if } 
\Hom(B,G) \overset{\sigma^*}{\arr} \Hom(N,G) \mbox{ is surjective }.$$
We claim that this is equivalent to say that
$$G \mbox{ is } G'-\mbox{complete.}$$

Suppose that $\sigma^*$ is surjective. Given a sequence 
$c_m \in G$ $(m \in \omega)$, we want to find 
solutions $e_m,d_i \in G$ such that
$$e_{m+1} p_m = e_m + \sum\limits_{i < n} d_i k_{i m} + c_m.$$
Recall from Observation \ref{3.1} (c) and the note in the
proof of (c) that
$$N = \left < y_{m+1}p_m - y_m - \sum\limits_{i < n} x_i k_{i m}: \ \
m \in \omega  \right >_R = 
\bigoplus\limits_{m\in \omega}(y_{m+1}p_m - y_m - \sum\limits_{i < n} x_i 
k_{i m})R\subseteq B.$$ 
We define a homomorphism $\varphi$ on the free generators
sending $(y_{m+1}p_m - y_m - \sum\limits_{i < n} x_i k_{i m})$ to $c_m$ 
for all $m \in \omega$. 
From surjectivity of $\sigma^*$ we find a homomorphism
$\Phi: B \arr G$ which coincides with $\varphi$ on (the free generators of) 
$N$. If we put $\Phi(y_m) = e_m$ and $\Phi(x_i) = d_i$, then the last 
displayed equation holds and $G$ is $G'-$complete.

\bigskip
Conversely, let $G$ be $G'$-complete and let 
$$ 0 \arr G \overset{\gamma}{\arr} H \overset{\eta}{\arr} G' \arr 0$$
be a short exact sequence. We
want to construct the splitting map $\sigma: G' \arr H$ satisfying
$ \sigma \eta = i d_{G'}$.
Following Definition \ref{3.3} we define $\sigma$ first on the module freely 
generated by the $y'_m, {x'_i} 's$ and show that $N$ is mapped to 0, hence 
$\sigma$ will be well-defined on $G'$. 
Let $x^*_i, y^*_m \in H$ be preimages of $x'_i, y'_m \in G'$
under $\eta$ in the short exact sequence, hence 
$x^*_i \eta = x_i $  $(i < n)$ and $y^*_m \eta = y'_m$.
The relations in Definition \ref{3.3} viewed in $H$ give elements 
$c_m \in G$ $(m \in \omega)$ such that
$$y^*_{m+1} p_m = y^*_m + \sum\limits_{i < n} x_i^* k_{i m} + c_m
\gamma.$$
Recall that $G$ is $G'$-complete, hence we can find $e_m, d_i \in G \;\;(m \in \omega,
i < n)$ such that
$$e_{m+1} p_m = e_m + \sum\limits_{i < n} d_i k_{i m} + c_m.$$
Now we correct our first choice of elements and define
$$x_i'' \sigma' = x^*_i - d_i \gamma \quad (i < n), \quad \quad
y''_m \sigma' = y^*_m - e_m \gamma \quad (m \in \omega)$$
which is defined on the `canonical' free resolution $F$ of $G'$,
generated by elements $x_i'', y_m''$. 

We must show that the relations $N \subseteq F$ defining $G' = F/N$
are mapped to $0$. An arbitrary generator of $N$ is of the form
$$w = y''_{m+1} p_m - y''_m - \sum\limits_{i < n} x''_i k_{i m}.$$
We apply $\sigma$ and derive,
$$w \sigma' = (y^*_{m+1} - e_{m+1} \gamma) p_m - (y^*_m - e_m
\gamma)
- \sum\limits_{i<n} (x^*_i - d_i \gamma) k_{i m}$$
$$ = y^*_{m+1} p_m - (e_{m+1} p_m) \gamma - y^*_m + e_m \gamma - 
\sum\limits_{i < n} x^*_i k_{i m} + (\sum\limits_{i < n} d_i k_{i m}) 
\gamma$$
$$ = y^*_m + \sum\limits_{i<n}x^*_i k_{i m} + c_m \gamma - e_m \gamma - 
(\sum\limits_{i<n} d_i k_{i m})\; \gamma - c_m \gamma $$
$$- y^*_m + e_m \gamma - \sum\limits_{i<n} x^*_i k_{i m} + 
(\sum\limits_{i<n} d_i k_{i m})\gamma  = 0.$$
Hence $\sigma'$ induces $\sigma: G' \arr H$.
Obviously $\sigma \eta = id$ on the generators of $G'$, hence $\sigma \eta =
id_{G'}$ and $\Ext (G', G) = 0$.
\hfill ${\square}$

\section{Applications of the Ext-Lemma}

In this section we will give first applications of our Ext-Lemma \ref{3.4} which
are important later on. We will begin rederiving some known results due to
Harrison, Kaplansky and Salce, see \cite{Fu} and \cite{Sa}.

Let $\rho = (r_n)_{n \in \omega}$ be a divisibility chain of
positive integers $r_n$  ($r_0 = 1$ and $r_n | r_{n+1}$, 
say $r_{n+1} =q_n r_n$ for all $n \in \omega$).
If $G$ is a torsion-free abelian group,
then we define\\
$(i)$ \quad  $\Z_{(\rho)} = \left< 1/r_n : n \in \omega 
\right > \subseteq \Q$, the rational subgroup of $\Q$ generated by the $ 
1/r_n$'s.\\
$(ii)$ \quad  the $\rho$-topology on $G$ to be generated by the 
open sets $Gr_n$ for all $n \in \omega$\\
$(iii)$ \quad  $G_\rho = \bigcap\limits_{n \in \omega} G{r_n}$.

The sequence $\rho$ is essentially the characteristic of the
rational group $\Z_{(\rho)}$. If $\rho$ runs over all prime powers
different from a fixed prime $p$, then $\Z_{(\rho)} = \Z_{(p)}$ and
if $r_n = p^n$ for all $n \in \omega$, then $\Z_{(\rho)} = \Q^{(p)}$.
For obvious choices of $\rho$ in $(ii)$ we obtain the $p$-adic and
the $\Z$-adic topology, respectively. The $\rho$-topology on $G$ is
Hausdorff if and only if $G_\rho = 0$ in $(iii)$. We say that $G$ is
$\rho$-reduced. Note that $(G/G_\rho)_\rho = 0$, hence $G_\rho$ is a
radical and $G/G_\rho$ is $\rho$-reduced.

Recall that $G$ is $\rho$-complete if $G$ is complete in the
$\rho$-topology. We now apply our Ext-Lemma. 

\bigskip

\begin{corollary} \label {4.1} Let $G$ be torsion-free and
$\rho$-reduced for some $\rho$. Then

$$\Ext (\Z_{(\rho)},G) = 0 \iff G \mbox{ is complete in the }
\rho \mbox{-topology.}$$
\end{corollary}

{\bf Remark:} If $\rho = (n!)_{n \in \omega}$, then $\Z_{(\rho)} = \Q$
and the $\rho$-topology is the $\Z$-topology. Then (\ref{4.1}) is 
due to Harrison and Kaplansky, see \cite[p.235]{Fu}.
If $\rho = (p^n)_{n \in \omega}$, then the
$\rho$-topology is the $p$-adic topology and $G$ is $p$-adically
complete by (\ref{4.1}), see Salce \cite[Theorem 3.5]{Sa} for a
general discussion.

\bigskip

{\bf Proof.} Let $\rho$ be as above and put $y_m' = 1/r_n$, hence
$$ y_{m+1}' q_n = y'_m \mbox{ and }  q_m r_m = r_{m+1} $$
and  $\Z_{(\rho)}$ is a $0$-free-by-$1$ $\Z$-module by (3.2).
By the Ext-Lemma \ref{3.4} we note that $\Ext (\Z_{(\rho)},G) = 0$ 
is equivalent to say that $G$ is $\Z_{(\rho)}$-complete 
in the sense of Definition \ref{3.3}. 
Hence any sequence $c_n \in G$  gives rise to solution $y_n \in G$  
of the equations
$$y_{n+1}q_n = y_n + c_n$$
for ($n \in \omega$).
Using $(\rho)$ above we see that 
$y_0 = y_1r_1 -c_0 = y_2r_2 - (c_0 + c_1r_1)$
and inductively it follows that
$y_0 = y_nr_n - \sum^{n-1}\limits_{i=0} c_i r_i$, hence
$- y_0 = \sum\limits_{i \in \omega} c_i r_i \in G $. 
Any sequence $c_n \in G$ has a limit 
$\sum\limits_{i \in \omega} c_i r_i \in G $ and $G$ is complete in
the $\rho$-topology. The converse is obvious and the Corollary 4.1 is
shown.  $\hfill{\square}$

\bigskip

The following definition extends the notion of a splitter.

\bigskip

\begin{definition} \label{4.2} We say that the $R$-module $G$ over its 
nucleus $R$ is a finite-rank splitter if and only if $\Ext(G', G) = 0$ 
for all finite rank $R$-submodules $G'$ of $G$.
\end{definition}

Note that splitters are finite rank splitters. 
The existence of non-free but $\aleph_1$-free splitters mentioned 
in the introduction show that the converse does not hold; 
see also \cite{GS,GT}.

\begin{proposition} \label{4.3} If $U \neq 0$ is a pure subgroup 
of a finite-rank
splitter $G$ and $p^{-1} \in \nuc U \setminus \nuc G$, then 
$\Ext(\Q^{(p)}, G) = 0$.

\end{proposition}

{\bf Proof.} Let $\bar{G} = G/p^\omega G$ where 
$p^\omega G = \bigcap\limits_{n \in \omega}
G p^n$ and note that $p^\omega \bar{G}$ is $p$-divisible and $\bar{G}$ is
$p$-reduced. Let
$$G' = \langle y_m : m \in \omega, y_{m+1} p = y_m \rangle \subseteq U,$$
hence $G' \cong \Q^{(p)}$. 
If $G$ is a finite-rank splitter, then $\Ext(G', \bar{G}) = 0$
and by (\ref{4.1}) $\bar{G}$ 
is a complete module over the $p$-adic integers $J_p$. Hence
$\bar{G}$ is cotorsion, see Fuchs \cite[p. 163]{Fu}. Now it is easy to 
see that $\Ext(\Q^{(p)}, G) = 0$
as well; note that $\bar{G}$ is $q$-divisible by all 
primes $q \neq p$; e.g. apply
Salce \cite[p. 21]{Sa}, see Theorem \ref{2.1}. $\hfill{\square}$

\medskip

We have an immediate 

\begin{corollary} \label {4.4} If $G$ is a finite-rank splitter which 
is $p$-reduced for
all primes $p$ and $0 \neq U$ is pure in $G$, then $\nuc U = \nuc G$.
\end{corollary}

\bigskip

{\bf Proof.} If $p^{-1} \in \nuc U \setminus \nuc G$, 
then $\Ext (\Q^{(p)}, G) = 0$ by (\ref{4.3})
and $G$ is a $J_p$-module which contradicts that $G$ is $q$-reduced for
$q \neq p$. $\hfill{\square}$

\bigskip

\begin{remark} \label {4.5} Note that $G = J_p \oplus J_q  \ (p \neq q)$ 
is a splitter with $\nuc J_p \neq \nuc G.$ 
Hence the hypothesis in (\ref{4.4}) cannot be dropped.
\end{remark}

Corollary \ref{4.1} might support the conjecture that a similar completeness, e.g.
for a different topology, would follow for (non-free) 
$n$-free-by-$1$ groups in place
of $\Z_{(\rho)} \subseteq \Q$.

\medskip
The following Theorem \ref{4.6} however shows that such a conjecture fails 
dramatically. The same example and modifications will serve for a different
purpose later on as well.
We begin with a definition  which generalizes $G'$-complete modules 
in order to deduce a result on rational cotorsion theories in Section 6.

\begin{definition} \label{5.1}
Let $\Phi$ be a set of finite-rank-free-by-$1$ $R$-modules and
$\F$ the class of all free $R$-modules
over some proper subring $R \subset \Q$.\\
(a) The $R$-module $L$ is $\Phi$-complete if and only if $\Ext (G',L) = 0$
for all $G' \in \Phi$.\\
(b) An $R$-module $G$ is $\Phi$-represented if $G$ can be written as
$\bigcup\limits_{\alpha< \lambda} G_\alpha$ a union of an 
ascending continuous chain of $R$-submodules 
$G_\alpha$ with $G_0 = 0$ and $G_{\alpha +1}/G_\alpha \in \Phi \cup \F$ 
for any $ \alpha \in \lambda $.
\end{definition}

\begin{remark} \label{5.2}
{\rm If $G$ is a torsion-free $R$-module, then we can find $\Phi$
such that $G$ is $\Phi$-represented: We define inductively $\Phi$ and
$\{G_\alpha : \alpha < \alpha^*\}$. If $G/G_\alpha$ is $\aleph_1$-free, then
choose any countable extension $G_{\alpha+1}$ of $G_\alpha$ which is pure in 
$G$ such that $G_{\alpha+1}/G_\alpha$ is free over $R$. 
Otherwise, there is a non-free, torsion-free, pure submodule 
$G_{\alpha+1}/G_\alpha
\subseteq G/G_\alpha$ of minimal rank $n$ by Pontryagin's theorem.
From Observation \ref{3.1} 
we infer that $G_{\alpha+1} / G_\alpha$ is an $n$-free-by-$1$ $R$-module 
which we add to $\Phi$. In case of limit ordinals, we just
take unions to define the next member of the chain.}
\end{remark}

\medskip

\begin{theorem} \label {4.6} Let $\Phi$ be a non-empty set of 
finite-rank-free-by-$1$  $R$-modules for some proper subring $R$ 
of $\Q$ and let $\kappa > |\Phi |$ be some infinite, regular cardinal with 
$\kappa = \kappa^{\aleph_0}$.
Then there exists a torsion-free $\Phi$-complete, $\Phi$-represented 
$R$-module $G$ of rank $\kappa$. If all modules in $\Phi$ have rank at 
least $n+1$, then any $R$-submodule of rank $\leq n$ in $G$ is free. 
\end{theorem}

\bigskip

{\bf Remark} The last theorem is of particular interest if $\Phi$ is a 
singleton. If the module in $\Phi$ has not rank $1$, then $G$ in 
Theorem \ref{4.6} has many free, pure submodules of rank at least 1, which
shows that $G$ can not be complete in its $p$-adic topology or any of its 
natural generalizations.

\bigskip

{\bf Proof.} 
We begin with a set $G$ of cardinality $\kappa$ only used for enumerating
all $\omega$-tuples of elements in $G$ to ensure that the final module
$G$ has solutions to all required equations linked to $\Phi$. Alternatively
we can enumerate all  $\omega$-tuples of elements in $G_\alpha$ of each 
submodule (with repetitions) while doing the transfinite construction of $G$.

Let $G_0 = 0$ and $G_1 = \bigoplus\limits_{i \in \kappa} e_i R$ be a
free $R$-module of rank $\kappa$ and choose a {\bf set} $G \supset G_1$ 
with
$|G \setminus G_1| = \kappa$ 
from which we will pick an ascending continuous chain
$G_\alpha$ of submodules $(\alpha < \kappa)$ with 
$|G \setminus G_\alpha| = \kappa$.
We also choose an enumeration 
$\bar{c} ^\alpha = (c^\alpha_n : n \in \omega)$
of $\omega$-tuples of $G \ (\alpha < \kappa)$ 
such that each $\bar{c} \in G^\omega$
appears $\kappa$ times, i.e. 
$|\{ \alpha \in \kappa : \bar{c} ^\alpha = \bar{c} \} | = \kappa.$

\medskip
Similarly we choose an enumeration of $\Phi$ by $X_\alpha$ 
($\alpha \in \kappa$) with $\kappa$-repetitions.

If $G_\alpha$ is constructed, then we distinguish two cases for constructing
$G_{\alpha+1}.$

\bigskip

{\bf Case 1:} There is a $c^\alpha_n \in G \setminus G_\alpha$ 
for some $n \in \omega$.
We set $G_{\alpha + 1} = G_\alpha \oplus R$.

\bigskip

{\bf Case 2:} If $c^\alpha_n \in G_\alpha$ for all $n \in \omega$, 
then we apply Lemma \ref{3.4}. There is an extension
$G_{\alpha+1} \supset G_\alpha$ such that
$$
G_{\alpha+1} / G_\alpha = \langle y^\alpha_m + G_\alpha, 
x^\alpha_i + G_\alpha : i<n \rangle \cong X_\alpha \ \
\mbox{ with } \ y^\alpha_{m+1} p_m = y^\alpha_m + 
\sum\limits_{i<n} x^\alpha_i k_{im} + c^\alpha_m $$
 where  $p_m^\alpha, k_{im}^\alpha \in R \; \; (m \in \omega)$  come from 
(3.2) applied to $ X_\alpha$.                                                                                                                                                                                                                                                                                                                       
Hence 
$$G_{\alpha+1} = 
\langle G_\alpha, y^\alpha_m, x^\alpha_i : i<n, m \in \omega \rangle.$$
Finally $G = \bigcup\limits_{\alpha \in \kappa} G_\alpha$
is the union of the continuous chain $\{G_\alpha : \alpha \in
\kappa\}$ and $G$ is torsion-free
of rank $\kappa$. 

If $\bar{c} = (c_m) \in G^\omega$ 
then there is some $\beta < \kappa$
with $\bar{c} \in G_\beta^\omega$ because $cf \kappa > \aleph_0$. 
By enumeration
we also find $\alpha > \beta$ such that $\bar{c} = \bar{c} ^\alpha$; 
the desired
solution for (\ref{3.3}) is in $G_{\alpha + 1}$ by construction of 
$G_{\alpha + 1}$.
Hence $G$ is $G'$-complete and $\Ext (G', G) = 0$ by (\ref{3.4}). 

Next we assume that all modules in $\Phi$ have rank at least $n+1$.
It remains to show
that any pure submodule $F$ of rank $\leq n$ is free. If this is shown for
$rk F = k < n$, then let $F$ be of rank $k + 1$. We can choose 
$\beta \in \kappa$ minimal with $F \subseteq G_\beta$. If $\beta = 0$, then 
$F \subseteq_* G_0$ which is free and our claim holds. If $\beta > 0$, then
$\beta$ cannot be a limit ordinal, hence $\beta = \alpha + 1$. 

If $\alpha$ belongs to Case 1, 
then $F \subseteq_* G_\alpha \oplus R$ and $F \not\subseteq G_\alpha$.
We see that $F/(G_\alpha \cap F) \cong F + G_\alpha / G_\alpha \cong R$ and
$(G_\alpha \cap F) \oplus g R = F$ with $r k (G_\alpha \cap F) \leq k$.
The induction hypothesis for $G_\alpha \cap F \subseteq_* G_\alpha$ 
completes this case. 

\medskip

If $\alpha$ belongs to Case 2, we argue similarly: $F$ is a pure 
submodule of 
$$G_{\alpha +1} = \langle G_\alpha, y^\alpha_m, x^\alpha_i: i < n, 
m \in \omega \rangle$$ 
and $F \not\subseteq G_\alpha$. We see that
$$F/(G_\alpha \cap F) \cong (F + G_\alpha)/ G_\alpha \subseteq 
G_{\alpha + 1} / G_\alpha \cong G'.$$
Obviously $rk (F/(G_\alpha \cap F)) \leq k < n$ and any subgroup of rank
$< n$ of $G'$ is free, hence $F/(G_\alpha \cap F)$ is free and splits 
$(G_\alpha \cap F) \oplus F' = F$ where $rk F' \geq 1.$ Induction completes
this case as well. $\hfill{\square}$

\bigskip

Finally we modify the proof of Theorem \ref{4.6} to get

\bigskip

\begin{theorem} \label{4.7}
Let $A$ be a ring with free additive group $A^+$ of cardinality 
$|A| < \kappa$ for some regular cardinal $\kappa = \kappa^{\aleph_0}$  
and let $G'$ be an n-free-by-1 abelian group for some $n > 0$. 
Then there exists an abelian torsion-free
$G'$-complete group $G$ of rank $\kappa$ such that $\End G \cong A.$
\end{theorem}

\bigskip

{\bf Proof.} We adopt the construction in the proof of 
Theorem \ref{4.6} for $\Phi = \{G'\}$ adding
intermediate steps from constructions of abelian groups with prescribed
endomorphism rings, based on 
the stationary version of Shelah's black box, see Franzen, G\"obel \cite{FG}.

We enumerate the
traps $\tau_\alpha$ ($\alpha \in \lambda^*$) of the black box such that 
the construction of the module depends on the norm $|\tau_\alpha|$  
of the trap $\tau_\alpha $ which takes values in $\kappa$. 
Let $S_1, S_2$ be two disjoint stationary subsets of $\kappa$ and assume that
the enumeration of $\omega$-tuples $\bar{c} ^\alpha$ in (\ref{4.6}) 
uses only
$\alpha \in S_1$ as indexing set rather than $\kappa$. We assume the
reader to be familiar with the construction in \cite{CG} or in \cite{FG}. 
The modifications will be quite obvious.

\bigskip

Let $G_0 = \bigoplus\limits_{i \in \kappa} e_i A$ 
and note that $U = \widehat{G}_0$, the
$p$-adic completion of $G_0$, provides enough space to carry out the proof
given in (\ref{4.6}).  Moreover 
$|U^\omega| = |U|^{\aleph_0} = 
| G_0 |^{\aleph_0} = \kappa^{\aleph_0} = \kappa$,
and the mentioned enumeration of $\omega$-tuples is settled. We consider
$G'' = G' \otimes A$ 
which is the direct sum of $\rk A$ copies of $G'$. If $\beta$
is a limit ordinal and all $A$-modules $G_\alpha (\alpha < \beta)$ are 
constructed, then we take $G_\beta = \bigcup\limits_{\alpha < \beta} G_\alpha$.
If $\beta = \alpha + 1$, then we distinguish three main cases.

\bigskip

{\bf Case I.} If $|\tau_\alpha| \in \kappa \setminus (S_1 \cup S_2)$, then let
$$G_{\alpha + 1} = (G_\alpha \oplus e A)_* \subseteq U$$
for some suitable $e \in U$.

{\bf Case II.} If $|\tau_\alpha| \in S_1$, then either $\bar{c} ^\alpha \in G_\alpha^ \omega$
or not. If $c^\alpha_n \notin G_\alpha$ for some $n$, then we apply Case I.
Otherwise let $G^\bullet_{\alpha + 1}$ be an extension of $G_\alpha$ in $U$ such
that $G_{\alpha + 1} / G_\alpha \cong G''$ 
with ``solutions'' $x^\alpha_i, y^\alpha_m \in G_{\alpha + 1}$ of 
$$y^\alpha_{m+1} p_m = y^\alpha_m + 
\sum\limits_{i < n} x^\alpha_i k_{im} + c^\alpha_m \;\; (m \in \omega)$$
as in (\ref{4.6}). 
Then we take $G_{\alpha+1} = (G^\bullet_{\alpha+1})_* \subseteq U$.

{\bf Case III.} If $|\tau_\alpha| \in S_2$, then we follow \cite{CG} 
or \cite{FG}. Either the trap
 $\tau_\alpha$ is of no interest (the trap does not determine a partial
homomorphism or the partial homomorphism is scalar maltiplication by some
$a \in A$ ), then we apply Case I or $\alpha$ provides an unwanted partial
homomorphism $\varphi_\alpha$. 
In this case we let $G_{\alpha+1} = \langle G_\alpha, e A \rangle_*$
for some suitable $e \in U$ with $e \varphi_\alpha \notin G_{\alpha + 1}$. This
can be arranged in such a way that the support of $e$ is almost disjoint from
$G_\alpha$, see \cite{CG} for elements with branch-like support. Standard
arguments and the proof of (\ref{4.6}) ensure Theorem \ref{4.7}.

\bigskip

Combining (\ref{4.6}) and (\ref{4.7}) with Lemma \ref{3.4} we have the following

\begin{corollary} \label{4.8}
Let $G'$ be an n-free-by-1 group and $\kappa = \kappa^{\aleph_0}$ some
infinite, regular cardinal. Then the following holds.

(a) There exists a torsion-free abelian 
$p$-reduced group $G$ of rank $\kappa$ with
$\Ext (G', G) = 0$ which is not complete in its $p$-adic topology.

(b) If $A$ is any ring with free additive group $A^+$ of rank $< \kappa$,
then there exists a family of $2^\kappa$ non-isomorphic torsion-free
$A$-modules $G$ of rank $\kappa$ with $\End G = A$ and $\Ext (G', G) = 0.$

(c) There exists a family of $2^\kappa$ indecomposable, 
pairwise non-isomorphic abelian 
groups $G$ of rank $\kappa$ such that $\Ext (G', G) = 0.$

\end{corollary}

{\bf Proof.} (c) follows from (b) for $A = \Z$. (a) and (b) follow from the
preceding results.

\bigskip

For applications in Section 5 we note that a modification of the proof of
(\ref{4.7}) leads to

\begin{theorem} \label{4.9}
Let $A$ be a ring with free additive group $A^+$ of rank
$< \kappa$ for some regular cardinal $\kappa = \kappa^{\aleph_0}$ 
and let $G'$ be an n-free-by-1
group for some $n > 0$. Then there exists a torsion-free $G'$-complete $A$-module
$G = \bigcup\limits_{\alpha \in \kappa} G_\alpha$ with ascending, continuous
chain $\{G_\alpha : \alpha \in \kappa\}$ of pure $A$-submodules $G_\alpha$
such that $G_0 = \bigoplus\limits_{i \in \kappa} e_i A, \ G_{\alpha+1} / G_\alpha$
either isomorphic to $A$ or to $G' \otimes A = \bigoplus\limits_{rkA} G'$ and 
$\End\; G = A$.
\end{theorem}

\bigskip

{\bf Proof.} The only relevant change in the proof of (\ref{4.7}) is in Case III
(Case I is similar but simpler). Note that we can choose (besides $e = e_0 \in U$)
additional elements $e_1, \ldots, e_{n-1}$ $A$-independent modulo $G_\alpha$
such that 
$$e \varphi_\alpha \notin \langle G_\alpha, e_0  A, \ldots, e_{n-1} A; y_m A:m \in \omega \rangle =: G_{\alpha+1}$$
where 
$$y_{m+1} p_m \equiv y_m + \sum\limits_{i<n} e_i k_{im} \mod G_\alpha (m \in \omega).$$
Hence $G_{\alpha+1} /G_\alpha \cong G' \otimes A$, and we proceed as before.
$\hfill{\square}$

\bigskip

We close this section with a direct proof of (\ref{2.7}), 
as promised above. 

\bigskip

\begin{corollary} \label{4.10}
If $G$ is reduced and torsion-free of cardinality $< 2^{\aleph_0}$
but not $\aleph_1$-free, then $G$ is not a splitter.
\end{corollary}

{\bf Proof.} Suppose $G$ is a splitter, hence $\Ext (G,G) = 0$, and let
$R = \nuc G$. By Pontryagin's theorem there is an $R$-submodule $G'$ of $G$
of minimal finite rank, say $n$, which is not free. 
Let $G' = \langle y_m, \bigoplus\limits_{i<n} x_i R \rangle_R$
and note that 
$$y_{m+1} p_m = y_m + \sum\limits_{i<n} x_i k_{im} \ (m \in \omega)$$
holds as shown in (\ref{3.1}), (3.2). If $n \geq 1$, then $x_0$ exists 
and is pure in
$G'$. Clearly 
$G'/\bigoplus\limits_{i<n} x_i R$ has type $(p_m)_{m \in \omega}$
and is not $R$. 
In this case let $c_m = x_0$ and note that $p_m x = c_m$ 
has no solutions in $G'$. If $n = 0$, then $\nuc G = R$ forces $p_m G \neq G$
and the existence of $c_m \in G \; (m \in \omega)$ such that $p_m x = c_m$ has
no solution for all $m \in \omega$. We will use these elements \\
$(i)$ \quad $c_m \in G$ with no solution $x_m \in G$ for 
$p_m x_m = c_m \ (m \in \omega)$ \\
to construct a non-trivial element of 
$\Ext (G',G)$. If $v \in \ ^\omega \! \omega$, then let
$$F = \bigoplus\limits_{m \in \omega} y''_m R \oplus \bigoplus\limits_{i<n} 
x''_i R \oplus G$$ 
and 
$$N_v = \langle y''_{m+1} p_m - y''_m - \sum\limits_{i<n} 
x''_i k_{im} - c_m v (m) : m \in \omega \rangle_R.$$
If $H_v = F/N_v, \gamma:G \ra H_v\;\; (g \ra g + N_v)$, then $0 \arr G
\overset{\gamma}{\arr} H_v$
as in (3.3), and if $y''_m + N_v = y'_m, x'_i = x''_i + N_v$, then
$(x'_i \ra x_i), (y'_m \ra y_m)$ induces an endomorphism $\eta_v : H_v \ra G'$
such that \\
$(ii)$ \quad $0 \arr G \overset{\gamma}{\arr} H_v 
\overset{\eta_v}{\arr} G' \ra 0$
is exact.
 
As $\Ext (G', G) = 0$, $(ii)$ has a splitting map
$\sigma_v: G' \ra H_v$ such that $\sigma_v \eta_v = id_{G'}$. 

If $e^v_m = y'_m - y_m \sigma_v$ and $d^v_i = x'_i - x_i \sigma_v$ then $e^v_m, d^v_i \in \ker \eta_v$
and there are $e^v_m, d^v_i \in G$ with 
$e^v_m \gamma = e'_m, d^v_i \gamma = d'_i$.
The obvious relations in $G'$ and $H_v$ produce (by subtraction) new 
relations
$$e^v_{m+1} p_m = e^v_m + \sum\limits_{i<n} d_i^v k_{im} + c_m \gamma v(m), 
\ \ (v \in \;  ^\omega \omega).$$
Recall that $|G| < 2^{\aleph_0} = |^\omega \omega|$. By a pigeon-hole
argument there are $v \neq w \in \; ^\omega \omega$ with 
$d^v_i = d^w_i$
for all $i < n$; then $t$ is is defined to be the branch point of $v$ and 
$w$, and we may assume $t > n$. Note that
$v(m) = w (m)$ for all $m < t$ and $v(t) - w (t) =1$ without loss of
generality, moreover $e^v_m = e^w_m$ for all $m < t$. We subtract the
two sets of relations for $v$ and $w$ respectively and get 
$(e^v_{n+1} - e^w_{n+1}) p_t = c_t \gamma$. As $\gamma$ is a pure
embedding, the last equation contradicts our choice $(i)$ of $c_t$.

\section{Splitters which are neither free nor cotorsion}

In this section we want to answer Schultz's \cite[Problem 4]{Sch} in the 
negative by providing a list of splitters in ZFC which are neither free over
their nuclei nor cotorsion. Our Examples \ref{5.6} also show that there is no
hope of  classifying splitters because any prescribed 
$R$-algebra A which is free
as a $R$-module, $R$ a proper subring of $\Q$, is an endomorphism algebra 
$\End G \cong A$
of some splitter $G$ with nucleus $R$. Hence all kind of nasty decompositions
may occur, Kaplansky's test problems are violated etc., see Corner,
G\"obel \cite{CG}. We begin with a

\begin{remark} \label{5.3}
By the Ext-Lemma \ref{3.4} we see that $L$ in (\ref{5.1}) is $\Phi$-complete
if and only if all systems of equations $(G')$ related to $G' \in \Phi$ by
(3.2) have solutions in $G$.
\end{remark}

\begin{theorem} \label{5.4}
If $L$ is $\Phi$-complete and $G$ is $\Phi$-represented over 
the same nucleus, then $\Ext (G,L) = 0$
\end{theorem}

\begin{remark} We note that $L$ and $G$ must have the same nucleus $R$ which
follows by trivial modification of Corollary \ref{4.4}.
\end{remark}

{\bf Proof of 5.2}. Suppose $L$ and $G$ are $R$-modules as indicated and
$G = \bigcup\limits_{\alpha<\lambda} G_\alpha$ is the union of 
 an ascending, continuous chain
of $R$-modules over the nucleus $R$ with $G_{\alpha+1} / G_\alpha \in \Phi$
or $G_{\alpha + 1}/ G_\alpha$ free respectively.

\medskip

We must show that any sequence\\
$(i)$ \hfil $0 \arr L \arr H \overset{\sigma} \arr G \arr 0$\\
splits.

\medskip

We must find a splitting map $\eta: G \arr H$ with $\eta \sigma = id_G$. If
$H_\alpha = G_\alpha \sigma^{-1} \subseteq H$, then $H = \bigcup\limits_{\alpha<\lambda} H_\alpha$
and we construct $\eta$ by induction on $\alpha$, choosing an ascending continuous
chain $\eta_\alpha: G_\alpha \arr H_\alpha$ of splitting maps $\eta_\alpha \sigma_\alpha = id_{G_\alpha}$
for $\sigma_\alpha = \sigma \upharpoonright H_\alpha$. Hence $\eta = \bigcup\limits_{\alpha<\lambda} \eta_\alpha$
is as required. If $\eta_\alpha : G_\alpha \arr H_\alpha$ with $\eta_\alpha \sigma_\alpha = id_{G_\alpha}$
is given, we must find $\eta_{\alpha+1} \supset \eta_\alpha$ with $\eta_{\alpha+1}
\sigma_{\alpha+1} = id_{G_{\alpha+1}}$ for any $\alpha < \lambda$. If $G_{\alpha+1}/G_\alpha$
is free, then $G_{\alpha+1} = C_\alpha \oplus G_\alpha$ and $C_\alpha$ is a free
$R$-module. It is easy to define $\eta_{\alpha+1} \upharpoonright C_\alpha$ from
$\sigma_{\alpha+1}$, hence $\eta_{\alpha+1}$ is given component-wise. If
$G_{\alpha+1} / G_\alpha \in \Phi$, then we may assume that $G_{\alpha+1} /G_\alpha$
is given by (3.2).

\medskip
\noindent
$(ii)$ \quad $\langle y'_m, \bigoplus\limits_{i<n} x'_i R:y'_{m+1} p_m = y'_m +
\sum\limits_{i<n} x'_i k_{im}, m \in \omega\rangle_R$

\medskip
We take preimages $y_m, x_i \in G_{\alpha+1}$ of $y'_m, x'_i \;\;(m \in \omega, i<n)$
modulo $G_\alpha$. Note that $\sigma_{\alpha+1} : H_{\alpha+1} \arr G_{\alpha+1}$
is onto, hence we can take preimages $y''_m, x''_m \in H_{\alpha+1}$ under
$\sigma_{\alpha+1}$.

\medskip
\noindent
$(iii)$ $\begin{cases}
\ds y''_m \sigma = y''_m \sigma_{\alpha+1} = y_m 
& \mbox{ and } \quad x''_i \sigma = x''_i \sigma_{\alpha+1} = x_i\\
y_m + G_\alpha = y'_m 
& \mbox{ and } \quad x_i + G_\alpha = x'_i \quad \quad (m \in \omega, i < n)
\end{cases}$

\medskip
We want to find certain corrections $d_i \in L$ and $e_m \in L$ and
define a preliminary map $\eta_{\alpha+1}: G_{\alpha+1} \arr \widetilde{H}_{\alpha+1}$
where 
$H_{\alpha +1} \subseteq \widetilde{H}_{\alpha+1} = H_{\alpha+1} \otimes \Q$ denotes the divisible hull of 
$H_{\alpha+1}$. 
Note that $H$ is torsion-free as an extension of torsion-free
groups by $(i)$, hence $H_{\alpha+1}$ is torsion-free and the last
inclusion holds. We require $\eta_{\alpha+1} \upharpoonright G_\alpha = \eta_\alpha$ and
set\\
$(iv)$ \quad $x_i \eta_{\alpha+1} = x''_i + d_i$ and $y_m \eta_{\alpha+1} = y''_m + e_m
\quad (i < n, m \in \omega)$\\
The mapping $\eta_{\alpha+1}$ will be a well-defined homomorphism on
${G}_{\alpha+1}$ if the new relations are preserved. They come from
$(ii)$ and are of the form\\
$(v)$ \quad $y_{m+1} p_m = y_m + \sum\limits_{i<n} x_i k_{i m} + c^\alpha_m$ for some
$c^\alpha_m \in G_\alpha, m \in \omega.$\\
Our preliminary map $\eta_{\alpha+1}$ takes these equations to $H;$ by
$(iv)$ we can apply $\eta_{\alpha+1}$ to $(v)$ and get\\
$(vi)$ \quad $y''_{m+1} p_m + e_{m+1} p_m = 
y''_m + \sum\limits_{i<n} x''_i k_{i m}+e_m
+ \sum\limits_{i<n} d_i k_{i m} + c^\alpha_m \eta_\alpha$.\\
Hence $\eta_{\alpha+1}$ exists if we find solutions $e_m, d_i$ of
$(vi)$.\\
Put $c_m = -y''_{m+1} p_m + y''_m + \sum\limits_{i<n} x''_i k_{i m} +
c^\alpha_m \eta_\alpha$, note that $\eta_\alpha \sigma_\alpha = i d_{G_{\alpha}}$
and calculate with $(iii)$ and $(v)$
$$c_m \sigma = y_{m+1} p_m - y_m - \sum\limits_{i<n} x_i k_{i m} - c^\alpha_m.$$
Hence $c_m \in \ker \sigma = L$ for all $m \in \omega$ by $(v)$. 
The module $L$ is
$\Phi$-complete and in particular $\Ext (G_{\alpha+1} / G_\alpha, L) = 0$ holds.
By the Ext-Lemma \ref{3.4} we can find actual elements $d_i, e_m \in L$ such that
$$e_{m+1} p_m = e_m + \sum\limits_{i<n} d_i k_{i m} + c_m \quad \quad (m \in \omega).$$
Subtracting from $(vi)$ we obtain\\
$(vii)$ \quad $y''_{m+1} p_m = y''_m + \sum\limits_{i<n} x''_i k_{i m} - 
c_m + c^\alpha_m \eta_\alpha$\\
Our definition of $c_m$ above shows that $(vii)$ holds. We have seen that
$\eta_\alpha$ can be extended to $\eta_{\alpha+1}$ by $(iv)$ for suitable elements
$d_i, e_m$. Hence $\eta$ exists and is the derived splitting map of
$(i)$ by
$\eta \sigma = \bigcup\limits_\alpha \eta_\alpha \sigma_\alpha = \bigcup\limits_\alpha i d _{G_{\alpha}} = i d _G.$
\hfill$\square$

\bigskip

Theorem \ref{5.4} has two immediate consequences.

\bigskip

\begin{corollary} \label{5.5}
If $\Phi$ is a set of finite-rank-free-by-1 $R$-modules
and $G$ is a $\Phi$-complete, $\Phi$-represented $R$-module, then $G$ is a
splitter with nucleus $R$.
\end{corollary}

{\bf Proof}. Put $L = G$ and apply (\ref{5.4}).

\medskip                                 
If $\Phi = \{G'\}$ is a singleton and $G'$ is any n-free-by-1 abelian group
for some $n > 1$, then Theorem \ref{4.9} applies. Corollary \ref{5.5} provides a proper
class of nasty examples. 

\begin{example} \label{5.6}
If $\kappa = \kappa^{\aleph_0}$ is any cardinal and $n$ is any
integer $> 1$ and $A$ is a ring of cardinality $< \kappa$ with free additive
structure $A^+$, then there exists a splitter $G$ with $\End G \cong A$ such
that all its subgroups of rank $< n$ are free.
\end{example}

\bigskip

\begin{remark} \label{5.8}
(a) We also may prescribe the nucleus $R \subset \Q$
of the splitter $G$ in (\ref{5.6}) replacing the ring 
$A$ by an $R$-algebra $A$ with free
$R$-module structure $A_R$.\\
(b) Any example (\ref{5.6}) is a splitter but neither free nor cotorsion, hence
(\ref{5.6}) is a negative answer to the Problem 4 in Schultz \cite[p.11]{Sch}. The
examples exist in $ZFC$. 
\end{remark}

\bigskip
Problem 2 in Schultz
\cite{Sch} is the following question. If $G$ is a reduced torsion-free
splitter with nucleus $R$ such that $G$ has no countable homomorphic image
with nucleus $R$, is $G$ cotorsion? The construction of $G$ in Theorem \ref{4.9}
can be modified: Let $X$ be the direct sum of all countable torsion-free
reduced abelian groups with nucleus $\Z$ and note that $X$ is cotorsion-free
with $|X| = 2^{\aleph_0}$. In this case we do not prescribe the endomorphism
ring of $G$, but require $\Hom (G,X) = 0$. A by now standard argument \cite{CG}
and Shelah's black box show that $G$ exists satisfying all (other) conditions
in (\ref{4.9}). In this case $|G| = \kappa^{\aleph_0} = \kappa > 2^{\aleph_0}$, and
clearly $G$ serves as counter example for the above question. Note that
Problem 3 in \cite{Sch} was answered in G\"obel \cite{Go1}.

\section{Enough projectives and injectives in cotorsion theories generated
by sets of finite-rank-free-by-1 groups}

In this section we apply our methods concerning $G'$-complete groups $G$ from
Sections 4 and 5 for answering Problem 2 in Salce \cite[p.32]{Sa}. Recall the
definition of (rational) cotorsion theories and of the notion of ``enough
projectives'' from Section 2. Let $(\F, \C)$ be a cotorsion theory.
By Lemma \ref{2.2} it has enough projectives if any free groups $A$ gives rise to
$C \in \C, F \in \F$ and a short exact sequence

$$0 \arr A \arr C \arr F \arr 0.$$

\bigskip

A cotorsion theory $(\F, \C)$ is {\bf cogenerated by a set} if there
is a set of groups or equivalently a single group $X$ such that 
$\F = {}^\perp (X^\perp)$
and $\C = X^\perp$. Recall $(v)$, the notion ``cogenerated'' from Section 2!
Due to our knowledge of finite-rank-by-1 groups from Section 3 we will restrict
in this section to an arbitrary set $\Phi$ of finite-rank-by-1 groups (which
we may replace by their direct sum). We have a

\bigskip

\begin{mtheorem} \label{6.1}
If the cotorsion theory $(\F, \C)$ is cogenerated
by a set $\Phi$ of finite-rank-by-1 groups, then $(\F, \C)$ has
enough projectives and enough injectives.
\end{mtheorem}

\medskip

As a corollary (part (a)), we have the indicated answer of Salce's \cite{Sa}
problem. Moreover we are able to deal with ``quasi cotorsion'' and ``local
cotorsion''. Recall that a group $G$ is quasi cotorsion if 
$\Ext (\bigoplus\limits_p \Z_{(p)}, G) = 0$,
where $p$ runs over all primes. Dually
we define ``locally cotorsion'', see Salce \cite{Sa}. If $\Q^p$ is the 
subring 
$\{\frac{z}{p^{n}} : z, n \in \Z \}$ of $\Q$, then $G$ is locally
cotorsion if $\Ext (\bigoplus\limits_{p} \Q^p, G) = 0$. 
Obviously 
$\bigoplus\limits_{p} J_p \ (J_p = p$- adic integers) 
is locally cotorsion but
surely not cotorsion; more details are in \cite{Sa}. Note that $\Phi
= \{\Z_{(p)} : p$ any prime $\}$, or $\Phi = \{\Q^p$ : $p$ any prime$\}$, like
any rational group satisfies the hypothesis of (6.1).

\begin{corollary} \label{6.2} (a) All rational cotorsion theories have enough
projectives and enough injectives.\\
(b) The quasi cotorsion and the locally cotorsion theory have enough projectives
and enough injectives.
\end{corollary}

\bigskip

Remark: If the rational group in (a) is $\Q$, then Corollary (\ref{6.2}) (a)
is a classical result due to the founders of ``cotorsion'', 
see D. K. Harrison and details in Fuchs \cite{Fu}.

\bigskip

{\bf Proof of (\ref{6.1})}: By Salce's Lemma \ref{2.2} it is enough to begin
with a free abelian group $A$ and to construct $A' \supset A$ such that\\
$(i)$ \quad $\Ext (G', A') = 0$ for all $G' \in \Phi$ \\
and if $F = A'/A$ also\\
$(ii) \quad F \in \F_\Phi$ \\
or equivalently
$$F \in {}^\perp (\Phi)^\perp,  \mbox{ that is } \ (\Phi \perp X 
\Rightarrow F \perp X)$$
or more explicitly 
$$(\Ext(G', X) = 0 \quad \mbox{ for all } \quad G' \in \Phi) 
\Rightarrow \Ext (F,X) = 0.$$

We begin with $(i)$, the construction of $A'$. 
By Theorem \ref{4.6} we can find an extension 
$A' = \bigcup\limits_{\alpha < \lambda} A'_\alpha$ with 
$A'_0 = 0$ and $A'_1 = A$ such that $A'$ is $\Phi$-complete.

Hence $\Ext(G', A') = 0$ for all $G' \in \Phi$ and 
$A' \in \C_\Phi$ is shown. 

Now we show that 
$$0 \arr X \arr H \arr F \arr 0$$ 
splits.
Using the $\Ext$-Lemma \ref{3.4} and $\Ext(G',X) = 0$ for all
$G' \in \Phi$ we have that \\
$(iii)$ \quad  $ X$ is $G'$-complete for all $G' \in \Phi$. \\
Recall that $F = \bigcup\limits_{\alpha < \lambda}
F_\alpha$ with $F_\alpha = A'_\alpha/A'_1$ and $F_1 = 0$. Let $H_\alpha =
F_\alpha \sigma^{-1}$, hence $H = \bigcup\limits_{\alpha < \lambda} H_\alpha$
and $H_1 = X$ and we assume that $X \ra H$ in
the last exact sequence is the identity on $X$. 
The map $\sigma : H \ra F$ in this sequence gives
$\sigma_\alpha = \sigma \upharpoonright H_\alpha$
which induces\\
$(iv)$ \quad $H_{\alpha+1} / H_\alpha \cong F_{\alpha+1} / F_\alpha$\\
and inductively we will find splitting maps $\eta_\alpha$ such that \\
$(\alpha) \quad 0 \ra X \ra H_\alpha \overset{\sigma_\alpha}{\underset{\eta_\alpha}
{\rightleftarrows}} F_\alpha \ra 0$ and $\eta_\beta \subseteq \eta_\alpha$ 
for all $\beta < \alpha < \lambda.$

The desired splitting map will be $\eta = \bigcup\limits_{\alpha < \lambda} 
\eta_\alpha.$ If $\alpha = 1, F_1 = 0, \sigma_1 = 0$ and $\eta_1 = 0$ 
satisfies (1). If $\alpha$ is a limit and $(\beta)$ is defined for $\beta < \alpha$,
take $\eta_\alpha = \bigcup\limits_{\beta < \alpha} \eta_\alpha$ and $(\alpha)$
holds. It remains the case to construct $\eta_\alpha$
from $\eta_\beta$ for $\alpha = \beta + 1$ assuming $(\beta)$. If $F_\alpha/F_\beta$
is free, then $H_\alpha = H_\beta \oplus C$ for some free abelian group $C$. In
this case $\eta_\alpha$ is constructed easily from $\sigma_\alpha \upharpoonright C$.
Otherwise we may assume $F_\alpha/F_\beta \cong H_{\alpha+1} / H_\alpha \cong G'
\in \Phi$ from $(iv)$ and use the construction (\ref{4.6}) of $A'$. Let 
$$G' = \langle y'_m, x'_i: i < n, m \in \omega\rangle$$ 
with relations
$$y'_{m+1} p_m = y'_m + \sum\limits_{i<n} x'_i k_{im} \;\; (m \in \omega)$$
from (\ref{3.3}). Hence we can write
$$F_\alpha = \langle F_\beta, y''_m, x''_i : i < n, m \in \omega \rangle.$$
The above relations turn into\\
$(v)$\quad  $y''_{m+1} p_m = y''_m + \sum\limits_{i<n} x''_i k_{im} + c_m \;\;(m \in \omega)$
for some $c_m \in F_\beta$. 

The mapping $\sigma_\alpha$ is surjective by
$(\alpha)$ and we can find preimages $y_m, x_i \in H_\alpha$ and $e_m \in H_\beta$
such that $y_m \sigma_\alpha = y''_m, x_i \sigma_\alpha = x''_i$ and 
$e_m \sigma_\beta = c_m$. Note that 

$$H_\alpha = \langle H_\beta, y_m, x_i : i < n,  m \in \omega \rangle$$
by the isomorphism $(iv)$. Hence there are elements $e_m \in H_\beta$ such that
the above relations become \\
$(vi)$ \quad $y_{m+1} p_m = y_m + 
\sum\limits_{i<n} x_i k_{im} + e_m \;\;\; (m \in \omega)$.

\medskip

Like in the proof of Theorem 5.2 we note that 
$\eta_\alpha \supset \eta_\beta$
can be adjusted by elements in $X = \ker \sigma$. We claim that $\eta_\alpha$
can be achieved by taking\\
$(vii)$ \quad $y''_m \eta_\alpha = y_m + f_m$ and $x''_i \eta_\alpha = x_i + d_i \;\;\;\;
(i < n, m \in \omega)$ \\
for some particular elements $f_m, d_i \in X$. Then $\eta_\alpha$ is a
well-defined homomorphism if it preserves the relations $(v)$. If we apply
$\eta_\alpha$ given by $(vii)$ to $(v)$, we see that the relations are preserved
if\\
$(viii)$ \quad $y_{m+1} p_m + f_{m+1} p_m = y_m + \sum\limits_{i<n} x_i k_{im}
+ f_m + \sum\limits_{i<n} d_i k_{im} + c_m \eta_\beta$ \\
holds in $H_\alpha$.
Put $g_m = -y_{m+1} p_m + y_m + \sum\limits_{i<n} x_i k_{im} + c_m \eta_\beta$,
and note that $g_m \in \ker \sigma = X \;\; (m \in \omega)$ by the same 
calculation as in the proof of Theorem 5.2. By $(iii)$ there are elements
$f_m, d_i \in X$ with $f_{m+1} p_m = f_m + \sum\limits_{i<n} d_i k_{im} + g_m\;
(m \in \omega)$. Again as in (5.2), the new elements ensure that
$(viii)$ holds
and $\eta_\alpha$ exists.
\hfill$\square$

\section{Notes on splitters and on $\omega$-splitters of size $\aleph_1$}

 From
$\S 2$ we know that splitters $< 2^{\aleph_0}$ are $\aleph_1$-free, moreover
there are obvious splitters, the free $R$-modules and torsion-free cotorsion
groups - lets call them trivial splitters. There are non-trivial splitters of
size $2^{\aleph_0}$ which are not $\aleph_1$-free as shown in $\S 5$.
Hence it is natural to deal with $\aleph_1$-free splitters of cardinality
$\aleph_1$, a problem which became an independent topic, now separated in a
joint paper \cite{GS}. 
We mention from \cite{GS} that all $\aleph_1$-free splitter of size 
$\aleph_1$ are free.
Particular splitters are $\omega$-splitters as defined by Schultz
\cite{Sch}:

\begin{definition} \label{7.1} A group $G$ is an $\omega$-splitter if 
$\Ext (\bigoplus_\omega G, \bigoplus_\omega G) = 0.$ Recall that $\bigoplus_\omega G =
\bigoplus\limits_{n \in \omega} e_n G$ is a direct sum of $\omega$ copies of
$G$.
\end{definition}

First we apply the Ext-Lemma and rederive the following result due to 
Phil Schultz  \cite{Sch}.

\begin{proposition} \label{7.2} If $A$ and $G$ are two torsion-free 
abelian groups with the same nucleus such that
$\Ext (A, \bigoplus_\omega G) = 0$, then $A$ is $\aleph_1$-free.
\end{proposition}

\bigskip

{\bf Proof}. Note that $\nuc A = \nuc G = R \neq \Q$ by hypothesis of the proposition,
in particular $G \neq 0$. If $A$ is not $\aleph_1$-free, by Pontryagin's theorem
there exists $G'\subseteq_* A$ of minimal rank $n$ non-free, and $G'$ is an
$n$-free-by-1 $R$-module as considered in Section 3. From
$\Ext(A, \bigoplus_\omega G) = 0$ and ($*$) we have 
$\Ext(G'_i, \bigoplus_\omega G) = 0$. The $R$-module $G'$ gives rise to a 
representation (3.2) and equations\\
$(i)$ \quad $y_{m+1} p_m = y_m + 
\sum\limits_{i<n} x_i k_{im} \;\;\; (i<n, m \in \omega).$\\
The sequence $(p_m)_{m \in \omega}$ represents a non-trivial type of $G$
because of $\nuc A = \nuc G = R$ and (3.1). Hence 
$\bar{G} = G /\bigcap\limits_{m \in \omega} Gp_m \neq 0$ and 
$\Ext(G', \bigoplus_\omega \bar{G}) = 0$ and\\
$(ii)$ \quad $D = \bigoplus\limits_{n \in \omega} e_n \bar{G} = 
\bigoplus_\omega \bar{G}$ is $G'$-complete \\
by our Ext-Lemma 3.3. If $c_m \in D\; (m \in \omega)$ we have solutions
\\
$(iii)$ \quad $ y_m, x_i \in D \ \mbox{ such that } \  y_{m+1} p_m =
y_m + \sum\limits_
{i<n} x_i k_{im} + c_m.$

\medskip
In order to derive a contradiction, we choose some $0 \neq c  \in \bar{G}$
and note that $\bigcap\limits_{m \in \omega} \bar{G}p_m = 0$.
Then let $c_n = e_nc$ for all $n \in \omega$.
Choose
$n_0 \in \omega$ large enough such that $x_i \in D'=
\bigoplus\limits_{s<n_0}
e_s \bar{G}$ for any $i < n$ and modulo $D'$ the equations $(iii)$ become
$$y_{m+1} p_m \equiv y_m + c_m$$
Inductively we have $y_0 \equiv y_{n+1} q_n - (c_0 + \sum^n\limits_{i\geq1} c_i q_{i-1})$
where we use $q_i = \prod^i\limits_{i=0} p_i$. Note that $y_{n+1} q_n \ra 0\;\;
(n \ra \infty)$ in the Hausdorff topology induced by $ Dq_i$ on $D$, hence
$$- y_0 \equiv c_0 + \sum^\infty\limits_{i\geq1} c_i q_{i-1}$$
in the limit. The right hand side is obviously not in $D/D'$ by our choice of
the $c_i's$, while the left hand side is in the direct sum $D/D'$ because
$y_0 \in D$, a contradiction.
\hfill$\square$

\medskip
Next we use the main result of \cite{GS} mentioned and apply a simple 
but clever reduction due to Schultz \cite{Sch}.

\begin{lemma} \label{7.4} {\bf (\cite{Sch})}  Any torsion-free group with
nucleus $R$ has an epimorphic image of size 
$\leq 2^{\aleph_0}$ with the same nucleus.
\end{lemma}

{\bf Proof}. First we note that we may assume that the $R$-module $G$ in 
question is reduced. Let $S$ be the set of all primes $p$ with $Gp \neq G$.
If $p \in S$, let $\bar{G}^p$ be the $p$-adic completion of $G/\bigcap\limits_
{n \in \omega} Gp^n$ and $\eta_p : G \ra \bar{G}^p$ be the canonical projection
followed by the embedding into $\bar{G}^p$. Obviously there is a projection
$\pi_p$ of the $p$-adic module such that $\frac{1}{p}$ is not in the image of
$G \eta_p \pi_p$. If $\sigma = \prod\limits_{p \in S} \eta_p \pi_p$ then
$G \sigma \subseteq \prod\limits_{p \in S} \bar{G}^p$, and $G \sigma$ has the
same nucleus as $G$.

Combining these observations with some known fact we are able to extend
\cite{Sch} by showing the following

\begin{theorem} \label{7.5} (ZFC + $\Diam_\kappa$) All $\omega$-splitters
of cardinality $\leq \kappa$ are free over their nuclei.
\end{theorem}

Theorem \ref{7.5} follows from a more general

\begin{lemma} \label{7.6} (ZFC + $\Diam_\kappa$) If $A$ and $G$ are
torsion-free abelian groups with the same nucleus $R$ with $|A| \leq \kappa$
and if $\Ext(A, \bigoplus_\omega G) = 0$, then $A$ is $R$-free.
\end{lemma}

{\bf Proof}. First, using (\ref{7.4}), we replace $G$ above by a group of
cardinality
$\aleph_1$ with the same nucleus $R$. Hence \\
$(i)$ $\Ext(A, \bigoplus_\omega G) = 0$ and $|G| \leq \aleph_1$ \\
without loss of generality. The claim follows by induction on $\kappa = |A|$.
If $A$ is countable, then $A$ is free by (\ref{7.2}). Now we may assume that 
$\kappa \geq \aleph_1$ and Theorem 1.15 in Eklof, Mekler \cite[p.353]{EM}
applies for regular $\kappa$:

If $M = \bigoplus_\omega G$ and $\Ext(A,M) = 0$, then
$\Gamma_M A = 0$. Here we used the induction hypothesis that $A$ is $\kappa$-free,
which follows from ($*$). From $\Gamma_M A = 0$ we find a cub in $\kappa$ 
(which we may identify with $\kappa$) such that $A = 
\bigcup\limits_{\alpha<\kappa}A_\alpha$ is a $\kappa$-filtration and $0 = 
\Ext(A_{\alpha+1} / A_\alpha, M) =\Ext(A_{\alpha+1} / A_\alpha, 
\bigoplus_\omega G)$. Hence $A_{\alpha+1} / A_\alpha$
is free for all $\alpha \in \kappa$, again by induction hypothesis and $A$ is
$R$-free. If $\kappa$ is a singular cardinal, then we recall that $R$ is 
hereditary and Shelah's Singular-Compactness-Theorem applies; see e.g.
\cite{EM} [p.107, Theorem 3.5]. Hence $A$ is $R$-free in this case as well.
\hfill$\square$

\noindent
R\"udiger G\"obel \\
Fachbereich 6, Mathematik und Informatik \\
Universit\"at Essen, 45117 Essen, Germany \\
{\small e--mail: R.Goebel@Uni-Essen.De}\\
and \\ 
Saharon Shelah \\
Department of Mathematics\\ 
Hebrew University, Jerusalem, Israel \\
and Rutgers University, Newbrunswick, NJ, U.S.A \\
{\small e-mail: Shelah@math.huji.ae.il}

\end{document}